\documentclass[graybox]{svmult}
\usepackage{mathptmx}       
\usepackage{helvet}         
\usepackage{courier}        
\usepackage{makeidx}         
\usepackage{graphicx}        
\usepackage{multicol}        
\usepackage[bottom]{footmisc}
\usepackage{amsmath,amssymb,amsfonts}        

\usepackage{graphicx}
\usepackage{amsmath}
\usepackage{amssymb}
\usepackage{latexsym}
\usepackage{subfigure}
\usepackage{crop}
\usepackage{algorithmic}
\usepackage{algorithm}
\usepackage{multirow}
\usepackage{bm}
\usepackage{bbm}
\usepackage{enumerate}
\usepackage{framed} 
\usepackage[framed]{ntheorem}
\usepackage{url}
\usepackage[colorlinks = true, pdfstartview = FitV, linkcolor = blue, citecolor = blue, urlcolor = blue]{hyperref}
\usepackage{array}
\usepackage{paralist}
\usepackage{authblk}

\usepackage{cleveref}
\crefname{equation}{Eq.}{Eqs.}
\crefname{figure}{Fig.}{Figs.}
\usepackage{relsize}

\usepackage{fullpage}

\usepackage[sort,nocompress]{cite}

\usepackage{enumitem}
\setlist[enumerate,1]{leftmargin=*,wide=0em, noitemsep,nolistsep, label = {\bfseries \arabic*.}}
\setlist[itemize,1]{leftmargin=*,wide=0em, noitemsep,nolistsep}

\usepackage{titlesec}
\titleformat*{\section}{\large\bfseries}
\titleformat*{\subsection}{\large\bfseries}
\titleformat*{\subsubsection}{\large\bfseries}
\titleformat*{\paragraph}{\normalsize\bfseries}
\titleformat*{\subparagraph}{\normalsize\bfseries}

\newcommand {\uu}  { {\bf u} }

\newcommand {\qq}  { {\bf q} }

\newcommand {\mm}  { {\bf m} }

\newcommand {\ff}  { {\bf f} }

\newcommand {\dd}  { {\bf d} }

\renewcommand{\vec}[1]{\ensuremath{\mathbf{#1}}}

\newcommand{\grad}{\ensuremath {\vec \nabla}}

\newcommand{\defeq}{\mathrel{\mathop:}=}

\definecolor{forestgreen}{rgb}{0.13, 0.55, 0.13}

\newcounter{cmt}
\usepackage{tcolorbox} 
\tcbuselibrary{breakable}
\tcbuselibrary{skins}

\usepackage{listings} 

\definecolor{mygreen}{rgb}{0,0.6,0}
\definecolor{mygray}{rgb}{0.5,0.5,0.5}
\definecolor{mymauve}{rgb}{0.58,0,0.82}

\newcommand\modelspace{{\cal M}}

\makeindex             

\begin{document}
\title*{Optimization Methods for Inverse Problems}
\author{Nan Ye and Farbod Roosta-Khorasani and Tiangang Cui}
\institute{Nan Ye \at ACEMS \& Queensland University of TechnologyInstitute, \email{n.ye@qut.edu.au}
\and Farbod Roosta-Khorasani \at University of Queensland \email{fred.roostauq.edu.au}
\and Tiangang Cui \at Monash University \email{tiangang.cui@monash.edu}}
\maketitle{}

\abstract{
Optimization plays an important role in solving many inverse problems. Indeed, the task of inversion often either involves or is fully cast as a solution of an optimization problem. In this light, the mere non-linear, non-convex, and large-scale nature of many of these inversions gives rise to some very challenging optimization problems.
The inverse problem community has long been developing various techniques for
solving such optimization tasks. However, other, seemingly disjoint communities,
such as that of machine learning, have developed, almost in parallel, interesting alternative methods which might have stayed under the radar of the inverse problem community. In this survey, we aim to change that. In doing so, we first discuss current state-of-the-art optimization methods widely used in inverse problems. We then survey recent related advances in addressing similar challenges in problems faced by the machine learning community, and discuss their potential advantages for solving inverse problems. 
By highlighting the similarities among the optimization challenges faced by the inverse problem and the machine learning communities, we hope that this survey can serve as a bridge in bringing together these two communities and encourage cross fertilization of ideas.}

\section{Introduction}
\label{sec:introduction}

Inverse problems arise in many applications in science and engineering. The term ``inverse problem'' is generally understood as the problem of finding a specific physical property, or properties, of the medium under investigation, using indirect measurements. This is a highly important field of applied mathematics and scientific computing, as to a great extent, it forms the backbone of modern science and engineering. Examples of inverse problems can be found in various fields within medical imaging \cite{arridge1999optical,arridge1997optical,bertero2010introduction,rundell1997inverse,louis1992medical} and several areas of geophysics including mineral and oil exploration \cite{menke2012geophysical,aster2013parameter,bunks1995multiscale,russell1988introduction}. 

In general, an inverse problem aims at recovering the unknown underlying parameters
of a physical system which produces the available observations/measurements.
Such problems are generally ill-posed \cite{hadamard1902sur}.
This is often solved via two approaches: 
a Bayesian approach which computes a posterior distribution
of the models given prior knowledge and the data, or a regularized data
fitting approach which chooses an optimal model by minimizing an objective that
takes into account both fitness to data and prior knowledge.
The Bayesian approach can be used for a variety of downstream inference tasks, such as credible intervals for the
parameters; it is generally more computationally expensive than the data fitting approach. The computational attractiveness of data fitting comes at a cost: it can only produce a ``point'' estimate of the unknown parameters. However, in many applications, such a point estimate can be more than adequate.  

In this review, we focus on the data fitting approach.
Optimization algorithms are central in this approach as the recovery of the unknown parameters is formulated as an optimization problem.
While numerous works have been done on the subject, there are still many
challenges remaining, including scaling up to large-scale problems, dealing with
non-convexity.
Optimization constitutes a backbone of many machine learning applications~\cite{bottou2016optimization,domingos2012few}. Consequently, there are many related developments in optimization from the
machine learning community. However, thus far and rather independently, the machine learning and the inverse problems communities have largely developed their own sets of tools and algorithms to address their respective optimization challenges.
It only stands to reason that many of the recent advances by machine learning can be potentially applicable for addressing challenges in solving inverse problems.
We aim to bring out this connection and encourage permeation of ideas across these two communities.

In \Cref{sec:inv_prob}, we present general formulations for the inverse
problem, some typical inverse problems, and optimization algorithms commonly
used to solve the data fitting problem.
We discuss recent advances in optimization in \Cref{sec:advances}.
We then discuss areas in which cross-fertilization of optimization and inverse
problems can be beneficial in \Cref{sec:discussion}.
We conclude in \Cref{sec:conclusion}.

\section{Inverse Problems}
\label{sec:inv_prob}

An inverse problem can be seen as the reverse process of a forward problem,
which concerns with predicting the outcome of some measurements given a complete
description of a physical system.
Mathematically, a physical system is often specified using a set of model
parameters $\mm$ whose values completely characterize the system.
The model space $\modelspace$ is the set of possible values of $\mm$.
While $\mm$ is usually arise as a parameter function, in practice it is often
discretized as a parameter vector for the ease of computation, typically using
the finite element method, the finite volume method, or the finite difference
method.
The forward problem can be denoted as 
\begin{equation} \label{eq:forward_prob}
	\mm \to \dd = \ff(\mm),
\end{equation}
where $\dd$ are the error-free predictions, and the above notation is a
shorthand for 
$\dd 
= (\dd_1, \ldots, \dd_s) 
= (\ff_1(\mm), \ldots, \ff_s(\mm))$,
with $\dd_{i} \in \mathbb{R}^{l}$ being the $i$-th measurement.
The function $\ff$ represents the physical theory used for the prediction and is
called the forward operator.
The observed outcomes contain noises and relate to the system via the following
the observation equation
\begin{equation} \label{eq:obs}
	\dd = \ff(\mm) + \bm{\eta},
\end{equation}
where $\bm{\eta}$ are the noises occurred in the measurements.
The inverse problem aims to recover the model parameters $\mm$ from such noisy
measurements.

The inverse problem is almost always ill-posed, because the same measurements
can often be predicted by different models. 
There are two main approaches to deal with this issue.
The Bayesian approach assumes a prior distribution $P(\mm)$ on the model and a
conditional distribution $P(\bm{\eta} \mid \mm)$ on noise given the model.
The latter is equivalent to a conditional distribution $P(\dd \mid \mm)$ on
measurements given the model.
Given some measurements $\dd$, a posterior distribution 
$P(\mm \mid \dd)$ on the models is then computed using the Bayes rule
\begin{equation}
	P(\mm \mid \dd) 
		\propto P(\mm) P(\dd \mid \mm).
\end{equation}
Another approach sees the inverse problem as a data fitting problem that finds
an parameter vector $\mm$ that gives predictions $\ff(\mm)$ that best fit the
observed outcomes $\dd$ in some sense.
This is often cast as an optimization problem
\begin{equation} \label{eq:invopt}
	\min_{\mm \in \modelspace}\quad 
		\psi(\mm, \dd),
\end{equation}
where the misfit function $\psi$ measures how well the model $\mm$ fits the data
$\dd$.
When there is a probabilistic model of $\dd$ given $\mm$, a typical choice of
$\psi(\mm, \dd)$ is the negative log-likelihood. 
Regularization is often used to address the issue of multiple solutions, and 
additionally has the benefit of stabilizing the solution, that is, the
solution is less likely to change significantly in the presence of outliers
\cite{vogelbook,ehn1, archer1995some}. 
Regularization incorporats some \textit{a priori} information on $\mm$ in the
form of a regularizer $R(\mm)$ and solves the regularized optimization
problem
\begin{equation} \label{eq:reg_invopt1}
	\min_{\mm \in \modelspace}\quad 
		\psi_{R,\alpha}(\mm, \dd) 
			\defeq \psi(\mm, \dd) + \alpha R(\mm),
\end{equation}
where $\alpha > 0$ is a constant that controls the tradeoff between prior
knowledge and the fitness to data.
The regularizer $R(\mm)$ encodes a preference over the models, with preferred
models having smaller $R$ values.
The formulation in \Cref{eq:reg_invopt1} can often be given a \textit{maximum a
posteriori (MAP)} interpretation within the Bayesian framework~\cite{scharf1991statistical}. 
Implicit regularization also exists in which there is no explicit term $R(\mm)$
in the objective~\cite{hansen1998,doas, doas5,rieder2005,rieder2010,hanke1}.

The misfit function often has the form $\phi(\ff(\mm), \dd)$, which
measures the difference between the prediction $\ff(\mm)$ and the observation
$\dd$.
For example, $\phi$ may be chosen to be the Euclidean distance between
$\ff(\mm)$ and $\dd$.
In this case, the regularized problem takes the form
\begin{equation} \label{eq:reg_invopt2}
	\min_{\mm \in \modelspace}\quad 
		\phi_{R,\alpha}(\mm, \dd) 
			\defeq \phi(\ff(\mm), \dd) + \alpha R(\mm),
\end{equation}
This can also be equivalently formulated as choosing the most preferred model
satisfying constraints on its predictions
\begin{equation} \label{eq:reg_invopt3}
	\min_{\mm \in \modelspace}\quad R(\mm), \quad
	 \text{ s.t. }\quad \phi(\ff(\mm), \dd) \le \rho.
\end{equation}
The constant $\rho$ usually relates to noise and the maximum discrepancy between
the measured and the predicted data, and can be more intuitive than $\alpha$.

\subsection{PDE-Contrained Inverse Problems}
\label{sec:application}
For many inverse problems in science and engineering, the forward model is not
given explicitly via a forward operator $\ff(\mm)$, but often conveniently
specified via a set of partial differential equations (PDEs).
For such problems, \Cref{eq:reg_invopt2} has the form
\begin{equation} \label{eq:pde-invopt}
	\min_{\mm \in \modelspace, \uu}\quad 
		\phi(P \cdot \uu, \dd) + \alpha R(\mm), \quad
			\text{ s.t. }\quad c_{i}(\mm, \uu_{i}) = 0, \quad i=1, \ldots, s,
\end{equation}
where $P \cdot \uu 
= (P_1, \ldots, P_s) \cdot (\uu_1, \ldots, \uu_s)
= (P_1 \uu_1, \ldots, P_s \uu_s)$ with $\uu_{i}$ being the field in the $i$-th
experiment, $P_{i}$ being the projection operator that selects fields at
measurement locations in $\dd_{i}$ (that is, $P_{i} \uu_{i}$ are the predicted values
at locations measured in $\dd_{i}$), and $c_{i}(\mm, \uu_{i}) = 0$ corresponds to the
forward model in the $i$-th experiment.
In practice, the forward model can often be written as 
\begin{equation} \label{eq:lin-model} 
	{\cal L}_{i}(\mm) \uu_{i} = \qq_{i}, \quad i = 1, \ldots, s,
\end{equation}
where ${\cal L}_{i}(\mm)$ is a differential operator, and $\qq_{i}$ is a term that
incorporates source terms and boundary values.

The fields $\uu_1, \ldots, \uu_s$ in \Cref{eq:pde-invopt} and
\Cref{eq:lin-model} are generally functions in two or three dimensional spaces,
and finding closed-form solutions is usually not possible.
Instead, the PDE-constrained inverse problem is often solved numerically by
discretizing \Cref{eq:pde-invopt} and \Cref{eq:lin-model} using the finite
element method, the finite volume method, or the finite difference method.
Often the discretized PDE-constrained inverse problem takes the form
\begin{equation} \label{eq:discretized-constrained}
	\min_{\mm \in \modelspace, \uu}\quad 
		\phi(P \uu, \dd) + \alpha R(\mm), \quad
			\text{ s.t. }\quad L_{i}(\mm) \uu_{i} = \qq_{i}, \quad i=1, \ldots, s,
\end{equation}
where $P$ is a block-diagonal matrix consisting of diagonal blocks 
$P_{1}, \ldots, P_{s}$ representing the discretized projection operators,
$\uu$ is the concatenation of the vectors $\uu_1, \ldots, \uu_s$ representing
the discretized fields,
and each $L_{i}(\mm)$ is a square, non-singular matrix representing the
differential operator ${\cal L}_{i}(\mm)$.
Each $L_{i}(\mm)$ is typically large and sparse.
We abuse the notations $P$, $\uu$ to represent both functions and their
discretized versions, but the meanings of these notations will be clear from
context.

The constrained problem in \Cref{eq:discretized-constrained} can be written in
an unconstrained form by eliminating $\uu$ using $\uu_{i} = L_{i}^{-1} \qq_{i}$,
\begin{equation} \label{eq:discretized-unconstrained}
	\min_{\mm \in \modelspace}\quad 
		\phi(P L^{-1}(\mm) \qq, \dd) + \alpha R(\mm), 
\end{equation}
where $L$ is the block-diagonal matrix with $L_1, \ldots, L_s$ as the diagonal
blocks, and $\qq$ is the concatenation of $\qq_1, \ldots, \qq_s$. Note that, as in the case of~\eqref{eq:reg_invopt2}, here we have $ \ff(\mm) = P L^{-1}(\mm) \qq $.

Both the constrained and unconstrained formulations are used in practice.
The constrained formulation can be solved using the method of Lagrangian
multipliers.
This does not require explicitly solving the forward problem as in the
unconstrained formulation.
However, the problem size increases, and the problem becomes one of finding a
saddle point of the Lagrangian, instead of finding a minimum as in the
constrained formulation.

%

\subsection{Image Reconstruction}
\label{sec:image_reconstruction}
Image reconstruction studies the creation of 2-D and 3-D images from sets
of 1-D projections. 
The 1-D projections are generally line integrals of a function representing the
image to be reconstructed.
In the 2-D case, given an image function $f(x, y)$, the integral along the line
at a distance of $s$ away from the origin and having a normal which forms an
angle $\phi$ with the $x$-axis is given by the Randon transform
\begin{equation}
	p(s, \phi) 
		= \int_{-\infty}^{\infty} f(z \sin \phi + s \cos \phi,
			-z \cos \phi + s \sin \phi) dz.
\end{equation}

Reconstruction is often done via back projection, filtered back projection, or
iterative methods \cite{natterer2001mathematical,herman2009fundamentals}.
Back projection is the simplest but often results in a blurred reconstruction.
Filtered back projection (FBP) is the analytical inversion of the Radon transform and
generally yields reconstructions of much better quality than back projection.
However, FBP may be infeasible in the presence of discontinuities or noise.
Iterative methods generally takes the noise into account and assumes a
distribution on the noise.
The objective function is often chosen to be a regularized likelihood of the
observation, which is then iteratively optimized using the expectation
maximization (EM) algorithm.

\subsection{Objective Function}
\label{sec:obj}
One of the most commonly used objective function is the least squares criterion,
which uses a quadratic loss and a quadratic regularizer.
Assume that the noise for each experiment in~\eqref{eq:obs} is independently but
normally distributed, i.e.,
$\bm{\eta}_{i} \sim \mathcal{N}(0,\Sigma_{i} ), \forall i$, where $\Sigma_{i}
\in \mathbb{R}^{l \times l}$ is the covariance matrix.
Let $\Sigma$ be the block-diagonal matrix with $\Sigma_{1}, \ldots, \Sigma_{s}$
as the diagonal blocks.
The standard \textit{maximum likelihood} (ML) approach~\cite{scharf1991statistical}, leads
to minimizing the least squares (LS) misfit function 
\begin{equation} \label{eq:l2}
	\phi(\mm) \defeq \|\ff(\mm) - \dd\|_{\Sigma^{-1}}^2,
\end{equation}
where the norm $\|x\|_{A} = \sqrt{x^{\top} A x}$ is a generalization of the Euclidean
norm (assuming the matrix $A$ is positive definite, which is true in the case of
$\Sigma_{i}^{-1}$). In the above equation, we simply write the general misfit function
$\phi(\ff(\mm), \dd)$ as $\phi(\mm)$ by taking the measurements $\dd$ as fixed
and omitting it from the notation.
As previously discussed, we often minimize a regularized misfit function
\begin{equation} \label{eq:reg_l2}
	\phi_{R,\alpha}(\mm) \defeq  \phi(\mm) + \alpha R(\mm).
\end{equation}
The prior $R(\mm)$ is often chosen as a Gaussian regularizer
$R(\mm) 
= (\mm - \mm_{\text{prior}})^{\top} \Sigma_{m}^{-1} (\mm - \mm_{\text{prior}})$.
We can also write the above optimization problem as minimizing $R(\mm)$ under
the constraints
\begin{equation} \label{eq:objective_eq}
	\sum_{i=1}^{s} \|\ff_{i}(\mm) - \dd_{i}\| \leq \rho.
\end{equation}

The least-squares criterion belongs to the class of $\ell_p$-norm
criteria, which contain two other commonly used criteria:
the least-absolute-values criterion and the minimax criterion
\cite{tarantola2005inverse}.
These correspond to the use of the $\ell_{1}$-norm and the $\ell_{\infty}$-norm
for the misfit function, while the least squares criterion uses the
$\ell_{2}$-norm.
Specifically, the least-absolute-values criterion takes 
$\phi(\mm) \defeq \|\ff(\mm) - \dd\|_1$, 
and the minimax criterion takes
$\phi(\mm) \defeq \|\ff(\mm) - \dd\|_{\infty}$.
More generally, each coordinate in the difference may be weighted.
The $\ell_{1}$ solution is more robust (that is, less sensitive to outliers) than the
$\ell_{2}$ solution, which is in turn more robust than the $\ell_{\infty}$ solution
\cite{claerbout1973robust}.
The $\ell_{\infty}$ norm is desirable when outliers are uncommon but the data
are corrupted by uniform noise such as the quantization errors
\cite{clason2012fitting}.

Besides the $\ell_{2}$ regularizer discussed above, the $\ell_{1}$-norm is
often used too.
The $\ell_{1}$ regularizer induces sparsity in the model parameters, that is,
heavier $\ell_{1}$ regularization leads to fewer non-zero model parameters.

\subsection{Optimization Algorithms}
\label{sec:optimization}

Various optimization techniques can be used to solve the regularized data
fitting problem.
We focus on iterative algorithms for nonlinear optimization below as the
objective functions are generally nonlinear.
In some cases, the optimization problem can be transformed to a linear program.
For example, linear programming can be used to solve the least-absolute-values
criterion or the minimax criterion.
However, linear programming are considered to have no advantage over
gradient-based methods (see Section 4.4.2 in \cite{tarantola2005inverse}), and
thus we do not discuss such methods here.
Nevertheless, there are still many optimization algorithms that can be covered
here, and we refer the readers to
\cite{bjorck1996numerical,nocedal2006numerical}.

For simplicity of presentation, we consider the problem of minimizing a
function $g(\mm)$.
We consider iterative algorithms which start with an iterate $\mm_0$, and
compute new iterates using
\begin{equation}
	\mm_{k+1} = \mm_{k} + \lambda_{k} p_{k},
\end{equation}
where $p_{k}$ is a search direction, and $\lambda_{k}$ a step size.
Unless otherwise stated, we focus on unconstrained optimization.
These algorithms can be used to directly solve the inverse problem in
\Cref{eq:reg_invopt1}. 
We only present a selected subset of the algorithms available and have to omit
many other interesting algorithms. 

\medskip\noindent{\bf Newton-type methods.}
The classical Newton's method starts with an initial iterate $\mm_{0}$, and computes new
iterates using 
\begin{equation}
	\mm_{k+1} = \mm_{k} - \left(\grad^2 g(\mm_{k})\right)^{-1} \grad g(\mm_{k}),
\end{equation}
that is, the search direction is 
$p_{k} = - \left(\grad^2 g(\mm_{k})\right)^{-1} \grad g(\mm_{k})$, and the step
length is $\lambda_{k} = 1$.
The basic Newton's method has quadratic local convergence rate at a small
neighborhood of a local minimum.
However, computing the search direction $p_{k}$ can be very expensive, and thus
many variants have been developed. In addition, in non-convex problems, classical Newton direction might not exist (if the Hessian matrix is not invertible) or it might not be an appropriate direction for descent (if the Hessian matrix is not positive definite). 

For non-linear least squares problems, where the objective function $g(\mm)$ is
a sum of squares of nonlinear functions, the Gauss-Newton (GN) method is often
used~\cite{sun2006optimization}.
Extensions to more general objective functions as in \Cref{eq:l2} with covariance matrix $ \Sigma $ and arbitrary regularization as in \Cref{eq:reg_l2} is considered in~\cite{roszas}.
Without loss of generality, assume 
$g(\mm) = \sum_{i=1}^{s} (\ff_{i}(\mm) - \dd_{i})^{2}$.
At iteration $k$, the GN search direction $p_{k}$ is given by
\begin{equation} 
	\left( \sum_{i=1}^{s} {J}_{i}^{\top} {J}_{i} \right) p_{k}  = -\grad g,
\end{equation}
where the sensitivity matrix $J_i$ and the gradient $\grad g$ are given by
\begin{align}
	J_i &= \frac{\partial \ff_{i}}{\partial \mm}(\mm_{k}), \quad i = 1, \ldots , s,\\
	\grad g &= 2 \sum_{i=1}^s J_{i}^T(\ff_{i}(\mm_{k}) - \dd_{i}),
\end{align}
The Gauss-Newton method can be seen as an approximation of the basic Newton's
method obtained by replacing $\grad^2 g$ by $\sum_{i=1}^{s} J_i^{\top} J_{i}$.
The step length $\lambda_{k} \in [0, 1]$ can be determined by a weak line
search~\cite{nocedal2006numerical} (using, say, the Armijo algorithm starting
with $\lambda_{k} = 1$) ensuring sufficient decrease in $g(\mm_{k+1})$ as
compared to $g(\mm_{k})$.

Often several nontrivial modifications are required to adapt this prototype method
for different applications, e.g., dynamic regularization~\cite{doas1,hanke1,rieder2005,rieder2010} and more general~\emph{stabilized GN} studied~\cite{rodoas1,doas12}.
This method replaces the solution of the linear systems defining $p_{k}$ by $r$
preconditioned conjugate gradient (PCG) inner iterations, which costs $2r$
solutions of the forward problem per iteration, for a moderate integer value $r$. 
Thus, if $K$ outer iterations are required to obtain an acceptable solution then
the total work estimate (in terms of the number of PDE solves) is approximated
{\em from below} by $2(r+1) K s$. 

Though Gauss-Newton is arguable the method of choice within the inverse problem community, other Newton-type methods exist which have been designed to suitably deal with the non-convex nature of the underlying optimization problem include Trust Region~\cite{conn2000trust,xu2017newton} and the Cubic Regularization~\cite{xu2017newton,cartis2012evaluation}. These methods have recently found applications in machine learning~\cite{xu2017second}. Studying the advantages/disadvantages of these non-convex methods for solving inverse problems can be indeed a useful undertaking.

\medskip\noindent{\bf Quasi-Newton methods.}
An alternative method to the above Newton-type methods is the quasi-Newton
variants including the celebrated limited memory BFGS (L-BFGS)
\cite{liu1989limited,nocedal1980updating}.
BFGS iteration is closely related to conjugate gradient (CG) iteration. In particular, BFGS applied to a strongly convex 
quadratic objective, with exact line search as well as initial Hessian $ P $, is equivalent to preconditioned CG with preconditioner $ P $. However, as the objective function departs from being a simple quadratic, the number of iterations of L-BFGS could be significantly higher than that of GN or trust region. In addition, it has been shown that the performance of BFGS and its limited memory version is greatly negatively affected by the high degree if ill-conditioning present in such problems~\cite{romassn1, romassn2,pyrrm_ssn_nonuni}. These two factor are among the main reasons why BFGS (and L-BFGS) can be less effective compared with other Newton-type alternatives in many inversion applications~\cite{haber2004quasi}.

\medskip\noindent{\bf Krylov subspace method.}
A Krylov subspace method iteratively finds the optimal solution to an
optimization in a larger subspace by making use of the previous solution in a
smaller subspace.
One of the most commonly used Krylov subspace method is the conjugate gradient
(CG) method.
CG was originally designed to solve convex quadratic minimization problems of
the form $g(\mm) = \frac{1}{2} \mm^{\top} A \mm - b^{\top} \mm$.
Equivalently, this solves the positive definite linear system $A\mm = b$.
It computes a sequence of iterates $\mm_0, \mm_1, \ldots$ converging to the minimum
through the following two set of equations.
\begin{align}
	\mm_0 &=0, & r_0 &= b, & p_0 &= r_0, & \\
  \mm_{k+1} &= \mm_{k} + \frac{||r_k||_2^2}{p_{k}^{\top} A p_{k}} p_{k}, 
    &
  r_{k+1} &= r_{k} - \frac{||r_k||_2^2}{p_{k}^{\top} A p_{k}} A p_k,
    &
	p_{k+1} &= r_{k+1} + \frac{||r_{k+1}||_2^2}{||r_{k}||_2^2} p_{k}, & k \ge 0.
\end{align}
This can be used to solve the forward problem of the form 
$L_{i}(\mm) \uu_{i} = \qq_{i}$, provided that $L_{i}(\mm)$ is positive definite, which
is true in many cases.

CG can be used to solve the linear system for the basic Newton direction.
However, the Hessian is not necessarily positive definite and modification is
needed \cite{nocedal2006numerical}.

In general, CG can be generalized to minimize a nonlinear function $g(\mm)$ 
\cite{fletcher2013practical,dai2011nonlinear}.
It starts with an arbitrary $\mm_0$, and $p_1 = - \grad g(\mm_0)$, and computes
a sequence of iterates $\mm_1, \mm_2, \ldots$ using the equations below:
for $k \ge 0$,
\begin{align}
	\mm_{k+1} 
		&= \arg\min_{\mm \in \{\mm_k + \lambda p_k, \lambda \in \mathbb{R}\}} g(\mm), \\
	p_{k+1}
		&= -\grad g(\mm_{k+1}) + \beta_{k} p_k, 
		\qquad\text{ where }
		\beta_k = \frac{||\grad g(\mm_{k+1})||_2^2}{||\grad g(\mm_{k})||_2^2}.
\end{align}
The above formula for $\beta_k$ is known as the Fletcher-Reeves formula.
Other choices of $\beta_k$ exist.
The following two formula are known as the Polak-Ribiere and Hestenes-Stiefel
formula respectively.
\begin{align}
    \beta_{k}
		&= \frac{\langle \grad g(\mm_{k+1}) - \grad g(\mm_{k}),  \grad g(\mm_{k+1}) \rangle}
        {||\grad g(\mm_{k})||_2^2}, \\
    \beta_{k}
		&= \frac{\langle \grad g(\mm_{k+1}) - \grad g(\mm_{k}),  \grad g(\mm_{k+1}) \rangle}
        {p_{k}^{\top} (\grad g(\mm_{k+1}) - \grad g(\mm_{k}))}.
\end{align}
In practice, nonlinear CG does not seem to work well, and is mainly used
together with other methods, such as in the Newton CG method
\cite{nocedal2006numerical}.

\medskip\noindent{\bf Lagrangian method of multipliers.}
The above discussion focuses on unconstrained optimization algorithms, which are
suitable for unconstrained formulations of inverse problems, or unconstrained
auxiliary optimization problems in methods which solves the constrained
formulations directly.
The Lagrangian method of multipliers is often used to directly solve the
constrained version.
Algorithms have been developed to offset the heavier computational cost and
slow convergence rates of standard algorthms observed on the Lagrangian, which
is a larger problem than the constrained problem.
For example, such algorithm may reduce the problem to a smaller one, such as
working with the reduced Hessian of the Lagrangian~\cite{haber2000optimization},
or preconditioning~\cite{haber2001preconditioned,benzi2011preconditioning}.
These methods have shown some success in certain PDE-constrained optimization
problems.

Augmented Lagrangian methods have also been developed
(e.g.~\cite{ito1990augmented,abdoulaev2005optical}).
These methods constructs a series of penalized Lagrangians with vanishing
penalty, and finds an optimizer of the Lagrangian by successively optimizing the
penalized Lagrangians.

\subsection{Challenges}
\label{sec:challenges}

\noindent{\bf Scaling up to large problems.}
The discretized version of an inverse problem is usually of very large scale,
and working with fine resolution or discretized problems in high dimension is
still an active area of research.

Another challenge is to scale up to large number of measurements, which is
widely believed to be helpful for quality reconstruction of the model in
practice, with some theoretical support.
While recent technological advances makes many big datasets available, existing
algorithms cannot efficiently cope with such datasets.
Examples of such problems include electromagnetic data inversion in mining
exploration \cite{na,dmr,haasol,olhash}, seismic data inversion in oil
exploration \cite{fichtner,hel,rnkkda}, diffuse optical tomography (DOT)
\cite{arridge1999optical,boas}, quantitative photo-acoustic tomography (QPAT)
\cite{gaooscher,yuan}, direct current (DC) resistivity
\cite{smvoz,pihakn,haheas,HaberChungHermann2010,doas12}, 
and electrical impedance tomography (EIT) \cite{bbp,cin,van2013lost}.

It has been suggested that many well-placed experiments yield practical advantage in order to obtain reconstructions of acceptable quality. 
For the special case where the measurement locations as well as the discretization matrices do not change from one experiment to another, various approximation techniques have been proposed to reduce the effective number of measurements, which in turn implies a smaller scale optimization problem, under the unifying category of ``simultaneous sources inversion'' \cite{rodoas1,roszas,roosta2015randomized,haber2014simultaneous,kumar2014GEOPemc}. Under certain circumstances, even if the $P_{i}$'s are different across experiments (but $ L_{i} $'s are fixed), there are methods to transform the existing data set 
into the one where all sources share the same receivers, \cite{rodoas2}. 

\medskip\noindent{\bf Dealing with non-convexity.}
Another major source of difficulty in solving many inverse problems, is the high-degree of non-linearity and non-convexity in~\eqref{eq:forward_prob}. This is most often encountered in problems involving PDE-constrained optimization where each $ \ff_{i} $ corresponds to the solution of a PDE. Even if the output of the PDE model itself, i.e., the ``right-hand side'', is linear in the sought-after parameter, the solution of the PDE, i.e., the forward problem, shows a great deal of non-linearity. This coupled with a great amount of non-convexity can have significant consequences in the quality of inversion and the obtained parameter. Indeed, in presence of non-convexity, the large-scale computational challenges are exacerbated, multiple folds over, by the difficulty of avoiding (possibly degenerate) \emph{saddle-points} as well as finding (at least) a \emph{local minimum}. 

\medskip\noindent{\bf Dealing with discontinuity.}
While the parameter function of the model is often smooth, the parameter
function can be discontinuous in some cases.
Such discontinuities arise very naturally as a result of the physical properties
of the underlying physical system, e.g., EIT and DC resistivity, and require
non-trivial modifications to optimization algorithms, e.g.,~\cite{rodoas1,
doas12}. 
Ignoring such discontinuities can lead to unsatisfactory recovery results
\cite{tali,doas12,doasleit2010}.
The level set method \cite{osse} is often used to model discontinuous parameter
function.
This reparametrizes the discontinuous parameter function as a differentiable
one, and thus enabling more stable optimization \cite{doasleit2010}.

\section{Recent Advances in Optimization}
\label{sec:advances}

Recent successes in using machine learning to deal with challenging perception
and natural language understanding problems have spurred many advances in the
study of optimization algorithms as optimization is a building block in machine
learning. 
These new developments include efficient methods for large-scale optimization, 
methods designed to handle non-convex problems, methods incorporating the structural
constraints, and finally the revival of second-order methods.
While these developments address a different set of applications in machine
learning, they address similar issues as encountered in inverse optimization and
could be useful.
We highlight some of the works below.
We keep the discussion brief because numerous works have been done behind these
developments and an indepth and comprehensive discussion is beyond the scope of
this review.
Our objective is thus to delineate the general trends and ideas, and provide
references for interested readers to dig on relevant topics.

\medskip\noindent{\bf Stochastic optimization.}
The development in large-scale optimization methods is driven by the
availability of many large datasets, which are made possible by the rapid
development and extensive use of IT technology. 
In machine learning, a model is generally built by optimizing a sum of misfit on
the examples.
This finite-sum structure naturally invites the application of stochastic optimization algorithms. This is mainly due to the fact that stochastic algorithms recover the sought-after models more efficiently by employing small batches of data in each iteration, as opposed to the whole data-set.
The most well-known stochastic gradient based algorithm is the stochastic gradient descent (SGD).
To minimize a finite-sum objective function 
\begin{align}
\label{eq:finite_sum}
g(\mm) = \frac{1}{n} \sum_{i=1}^{n} g_{i}(\mm),
\end{align} 
in the big data regime where $ n \gg 1 $, 
the vanilla SGD performs an update 
\begin{equation}
	\mm_{k+1} = \mm_{k} - \lambda_{k} \grad g_{i_{k}}(\mm_{k}),
\end{equation}
where $i_{k}$ is randomly sampled from $1, \ldots, n$.
As compared to gradient descent, SGD replaces the full gradient 
$\grad g(\mm)$ by a stochastic gradient $g_{i_{k}}(\mm_{k})$ with its
expectation being the full gradient.
The batch version of SGD constructs a stochastic gradient by taking the average
of several stochastic gradients.

SGD is inexpensive per iteration, but suffers from a slow rate of convergence.
For example, while full gradient descent achieves a linear convergence rate for
smooth strongly convex problems, SGD only converges at a sublinear rate.
The slow convergence rate can be partly accounted by the variance in the
stochastic gradient.
Recently, variance reduction techniques have been developed, e.g.
SVRG \cite{johnson2013accelerating} and SDCA \cite{shalev2013stochastic}.
Perhaps surprisingly, such variants can achieve linear convergence rates on
convex smooth problems as full gradient descent does, instead of sublinear rates
achieved by the vanilla SGD.
There are also a number of variants with no known linear rates but have fast
convergence rates for non-convex problems in practice, e.g.,
AdaGrad~\cite{duchi2011adaptive}, RMSProp~\cite{tijmen2012rmsprop},
ESGD~\cite{dauphin2015equilibrated}, Adam~\cite{kingma2014adam}, and
Adadelta~\cite{zeiler2012adadelta}.
Indeed, besides efficiency, stochastic optimization algorithms also seem to be able to
cope with the nonconvex objective functions well, and play a key role in the
revival of neural networks as deep learning~\cite{jin2017escape,ge2015escaping, levy2016power}.

\medskip\noindent{\bf Nonconvex optimization.}
There is also an increasing interest in non-convex optimization in the machine
learning community recently.
Nonconvex objectives not only naturally occur in deep learning, but also occur
in problems such as tensor decomposition, variable selection, low-rank matrix
completion, e.g. see \cite{ge2015escaping,mazumder2011sparsenet,jain2013low}
and references therein.

As discussed above, stochastic algorithms have been found to be capable of
effectively escaping local minima.
There are also a number of studies which adapt well-known acceleration
techniques for convex optimization to accelerate the convergence rates of both
stochastic and non-stochastic optimization algorithms for nonconvex problems,
e.g.,~\cite{li2015accelerated,allen2016variance,reddi2016stochastic,sutskever2013importance}.

\medskip\noindent{\bf Dealing with structural constraints.}
Many problems in machine learning come with complex structural constraints.
The Frank-Wolfe algorithm (a.k.a. conditional gradient)
\cite{frank1956algorithm} is an algorithm for optimizing over a convex domain.
It has gained a revived interest due to its ability to deal with many structural
constraints efficiently.
It requires solving a linear minimization problem over the feasible set, instead
of a quadratic program as in the case of proximal gradient algorithms or
projected gradient descent.
Domains suitable for the Frank-Wolfe algorithm include simplices,
$\ell_p$-balls, matrix nuclear norm ball, matrix operator norm ball
\cite{jaggi2013revisiting}.

The Frank-Wolfe algorithm belongs to the class of linear-optimization-based
algorithms \cite{lan2016conditional,lan2017conditional}.
These algorithms share with the Frank-Wolfe algorithm the characteristic of
requiring a first-order oracle for gradient computation and an oracle for
solving a linear optimization problem over the constraint set.

\medskip\noindent{\bf Second-order methods.}
The great appeal of the second-order methods lies mainly in the observed empirical performance as well as some very appealing theoretical properties. For example, it has been shown that stochastic Newton-type methods in general, and Gauss-Newton in particular, can not only be made scalable and have low per-iteration cost~\cite{rodoas1,rodoas2,roszas,haber2000optimization,haber2012effective,doas12}, but more importantly, and unlike first-order methods, are very \emph{resilient} to many adversarial effects such as \emph{ill-conditioning}~\cite{romassn1,romassn2,pyrrm_ssn_nonuni}. As a result, for moderately to very ill-conditioned problems, commonly found in scientific computing, while first-order methods make effectively no progress at all, second-order counterparts are not affected by the degree of ill-conditioning. 
A more subtle, yet potentially more severe draw-back in using first-order methods, is that their success is tightly intertwined with \emph{fine-tunning} (often many) \emph{hyper-parameters}, most importantly, the step-size~\cite{berahas2017investigation}. In fact, it is highly unlikely that many of these methods exhibit acceptable performance on first try, and it often takes many trials and errors before one can see reasonable results. In contrast, second-order optimization algorithms involve much less parameter tuning and are less sensitive to the choice of hyper-parameters~\cite{berahas2017investigation, xu2017second}. 

Since for the finite-sum problem~\eqref{eq:finite_sum} with $ n \gg 1 $, the operations with the Hessian/gradient constitute major computational bottlenecks, a rather more recent line of research is to construct the inexact Hessian information using the application of \emph{randomized methods}. Specifically, for convex optimization, the stochastic approximation of the full Hessian matrix in the classical Newton's method has been recently considered in~\cite{byrd2011use, byrd2012sample, wang2015subsampled,pilanci2015newton, erdogdu2015convergence, romassn1, romassn2, pyrrm_ssn_nonuni, Agarwal2016SecondOS, mutny2016stochastic, ye2016revisiting, bollapragada2016exact, mutny2017parallel, berahas2017investigation,eisen2017large}. In addition to inexact Hessian, a few of these methods study the fully stochastic case in which the gradient is also approximated, e.g.,~\cite{romassn1, romassn2,bollapragada2016exact}. For non-convex problems, however, the literature on methods that employ randomized Hessian approximation is significantly less developed than that of convex problems. A few recent examples include the stochastic trust region~\cite{xu2017newton}, stochastic cubic regularization~\cite{xu2017newton,tripuraneni2017stochastic}, and noisy negative curvature method~\cite{liu2017noisy}. Empirical performance of many of these methods for some non-convex machine learning applications has been considered in~\cite{xu2017second}.

\section{Discussion}
\label{sec:discussion}

Optimization is not only used in the data fitting approach to inverse problems,
but also used in the Bayesian approach.
An important problem in the Bayesian approach is the choice of the
parameters for the prior.
While these were often chosen in a somewhat ad hoc way, there are studies which
use sampling \cite{fox2016fast,agapiou2014analysis},
hierarchical prior models \cite{calvetti2007gaussian,calvetti2008hypermodels}, 
and optimization \cite{bardsley2010hierarchical,liu2017approximate} methods to
choose the parameters.
While choosing the prior parameters through optimization has found some success,
such optimization is hard and it remains a challenge to develop effective
algorithms to solve these problems.

For inverse problems with large number of measurements, solving each forward
problem can be expensive, and the mere evaluation of the misfit function may
become computationally prohibitive. 
Stochastic optimization algorithms might be beneficial in this case, because the
objective function is often a sum of misfits over different measurements.

The data fitting problem is generally non-convex and thus optimization
algorithms may be trapped in a local optimum.
Stochastic optimization algorithms also provide a means to escape the local
optima. 
Recent results in nonconvex optimization, such as those on accelerated methods,
may provide more efficient alternatives to solve the data fitting problem.

While box constraints are often used in inverse problems because they are easier to
deal with, simplex constraint can be beneficial. 
The Frank-Wolfe algorithm provides a efficient way to deal with the simplex
constraint, and can be a useful tool to add on to the toolbox of an inverse
problem researcher.

\section{Conclusion}
\label{sec:conclusion}

The state-of-the-art optimization methods used in the inverse problem community
do not yet cope with some challenges, such as large-scale problems,
nonconvexity, very well.
At the same time, many progresses in optimization have been made in the machine
learning community.
While our discussion on the connections are nevertheless brief, we hope it has
pointed out some interesting avenues that the two communities can learn from
each other and work together.

\bibliographystyle{plain}
\bibliography{ref_Nan,ref_Fred}

\begin{thebibliography}{100}

\bibitem{abdoulaev2005optical}
Gassan~S Abdoulaev, Kui Ren, and Andreas~H Hielscher.
\newblock {Optical tomography as a PDE-constrained optimization problem}.
\newblock {\em Inverse Problems}, 21(5):1507 -- 1530, 2005.

\bibitem{agapiou2014analysis}
Sergios Agapiou, Johnathan~M Bardsley, Omiros Papaspiliopoulos, and Andrew~M
  Stuart.
\newblock Analysis of the gibbs sampler for hierarchical inverse problems.
\newblock {\em SIAM/ASA Journal on Uncertainty Quantification}, 2(1):511--544,
  2014.

\bibitem{Agarwal2016SecondOS}
Naman Agarwal, Brian Bullins, and Elad Hazan.
\newblock Second order stochastic optimization in linear time.
\newblock {\em arXiv preprint arXiv:1602.03943}, 2016.

\bibitem{allen2016variance}
Zeyuan Allen-Zhu and Elad Hazan.
\newblock Variance reduction for faster non-convex optimization.
\newblock {\em arXiv preprint arXiv:1603.05643}, 2016.

\bibitem{archer1995some}
G.~Archer and DM. Titterington.
\newblock On some bayesian/regularization methods for image restoration.
\newblock {\em Image Processing, IEEE Transactions on}, 4(7):989--995, 1995.

\bibitem{arridge1999optical}
S.~R. Arridge.
\newblock Optical tomography in medical imaging.
\newblock {\em Inverse problems}, 15(2):R41, 1999.

\bibitem{arridge1997optical}
S.~R. Arridge and J.~C. Hebden.
\newblock Optical imaging in medicine: Ii. modelling and reconstruction.
\newblock {\em Physics in Medicine and Biology}, 42(5):841, 1997.

\bibitem{aster2013parameter}
R.~C. Aster, B.~Borchers, and C.~H. Thurber.
\newblock {\em Parameter estimation and inverse problems}.
\newblock Academic Press, 2013.

\bibitem{bardsley2010hierarchical}
Johnathan~M Bardsley, Daniela Calvetti, and Erkki Somersalo.
\newblock Hierarchical regularization for edge-preserving reconstruction of pet
  images.
\newblock {\em Inverse Problems}, 26(3):035010, 2010.

\bibitem{benzi2011preconditioning}
Michele Benzi, Eldad Haber, and Lauren Taralli.
\newblock A preconditioning technique for a class of pde-constrained
  optimization problems.
\newblock {\em Advances in Computational Mathematics}, 35(2):149--173, 2011.

\bibitem{berahas2017investigation}
Albert~S Berahas, Raghu Bollapragada, and Jorge Nocedal.
\newblock {An Investigation of Newton-Sketch and Subsampled Newton Methods}.
\newblock {\em arXiv preprint arXiv:1705.06211}, 2017.

\bibitem{bertero2010introduction}
M.~Bertero and P.~Boccacci.
\newblock {\em Introduction to inverse problems in imaging}.
\newblock CRC press, 2010.

\bibitem{bjorck1996numerical}
{\AA}ke Bj{\"o}rck.
\newblock {\em Numerical methods for least squares problems}.
\newblock SIAM, 1996.

\bibitem{boas}
D.A. Boas, D.H. Brooks, E.L. Miller, C.~A. DiMarzio, M.~Kilmer, R.J. Gaudette,
  and Q.~Zhang.
\newblock Imaging the body with diffuse optical tomography.
\newblock {\em Signal Processing Magazine, IEEE}, 18(6):57--75, 2001.

\bibitem{bollapragada2016exact}
Raghu Bollapragada, Richard Byrd, and Jorge Nocedal.
\newblock Exact and inexact subsampled {N}ewton methods for optimization.
\newblock {\em arXiv preprint arXiv:1609.08502}, 2016.

\bibitem{bbp}
L.~Borcea, J.~G. Berryman, and G.~C. Papanicolaou.
\newblock High-contrast impedance tomography.
\newblock {\em Inverse Problems}, 12:835--858, 1996.

\bibitem{bottou2016optimization}
L{\'e}on Bottou, Frank~E Curtis, and Jorge Nocedal.
\newblock Optimization methods for large-scale machine learning.
\newblock {\em arXiv preprint arXiv:1606.04838}, 2016.

\bibitem{bunks1995multiscale}
C.~Bunks, F.~M. Saleck, S.~Zaleski, and G.~Chavent.
\newblock Multiscale seismic waveform inversion.
\newblock {\em Geophysics}, 60(5):1457--1473, 1995.

\bibitem{byrd2011use}
Richard~H. Byrd, Gillian~M. Chin, Will Neveitt, and Jorge Nocedal.
\newblock On the use of stochastic {H}essian information in optimization
  methods for machine learning.
\newblock {\em SIAM Journal on Optimization}, 21(3):977--995, 2011.

\bibitem{byrd2012sample}
Richard~H. Byrd, Gillian~M. Chin, Jorge Nocedal, and Yuchen Wu.
\newblock Sample size selection in optimization methods for machine learning.
\newblock {\em Mathematical programming}, 134(1):127--155, 2012.

\bibitem{calvetti2007gaussian}
Daniela Calvetti and Erkki Somersalo.
\newblock A gaussian hypermodel to recover blocky objects.
\newblock {\em Inverse problems}, 23(2):733, 2007.

\bibitem{calvetti2008hypermodels}
Daniela Calvetti and Erkki Somersalo.
\newblock Hypermodels in the bayesian imaging framework.
\newblock {\em Inverse Problems}, 24(3):034013, 2008.

\bibitem{cartis2012evaluation}
Coralia Cartis, Nicholas~IM Gould, and Philippe~L Toint.
\newblock {Evaluation complexity of adaptive cubic regularization methods for
  convex unconstrained optimization}.
\newblock {\em Optimization Methods and Software}, 27(2):197--219, 2012.

\bibitem{cin}
M.~Cheney, D.~Isaacson, and J.~C. Newell.
\newblock Electrical impedance tomography.
\newblock {\em SIAM Review}, 41:85--101, 1999.

\bibitem{claerbout1973robust}
Jon~F Claerbout and Francis Muir.
\newblock Robust modeling with erratic data.
\newblock {\em Geophysics}, 38(5):826--844, 1973.

\bibitem{clason2012fitting}
Christian Clason.
\newblock L$\infty$ fitting for inverse problems with uniform noise.
\newblock {\em Inverse Problems}, 28(10):104007, 2012.

\bibitem{conn2000trust}
Andrew~R Conn, Nicholas~IM Gould, and Ph~L Toint.
\newblock {\em Trust region methods}, volume~1.
\newblock SIAM, 2000.

\bibitem{dai2011nonlinear}
Y.~Dai.
\newblock Nonlinear conjugate gradient methods.
\newblock {\em Wiley Encyclopedia of Operations Research and Management
  Science}, 2011.

\bibitem{dauphin2015equilibrated}
Yann Dauphin, Harm de~Vries, and Yoshua Bengio.
\newblock Equilibrated adaptive learning rates for non-convex optimization.
\newblock In {\em Advances in Neural Information Processing Systems}, pages
  1504--1512, 2015.

\bibitem{doasleit2010}
K.~van~den Doel, U.~Ascher, and A.~Leitao.
\newblock Multiple level sets for piecewise constant surface reconstruction in
  highly ill-posed problems.
\newblock {\em Journal of Scientific Computation}, 43(1):44--66, 2010.

\bibitem{doas12}
Kees van~den Doel and Uri Ascher.
\newblock Adaptive and stochastic algorithms for {EIT} and {DC} resistivity
  problems with piecewise constant solutions and many measurements.
\newblock {\em SIAM J. Scient. Comput.}, 34:DOI: 10.1137/110826692, 2012.

\bibitem{domingos2012few}
Pedro Domingos.
\newblock A few useful things to know about machine learning.
\newblock {\em Communications of the ACM}, 55(10):78--87, 2012.

\bibitem{dmr}
O.~Dorn, E.~L. Miller, and C.~M. Rappaport.
\newblock A shape reconstruction method for electromagnetic tomography using
  adjoint fields and level sets.
\newblock {\em Inverse Problems}, 16, 2000.
\newblock 1119-1156.

\bibitem{duchi2011adaptive}
John Duchi, Elad Hazan, and Yoram Singer.
\newblock Adaptive subgradient methods for online learning and stochastic
  optimization.
\newblock {\em The Journal of Machine Learning Research}, 12:2121--2159, 2011.

\bibitem{eisen2017large}
Mark Eisen, Aryan Mokhtari, and Alejandro Ribeiro.
\newblock {Large Scale Empirical Risk Minimization via Truncated Adaptive
  Newton Method}.
\newblock {\em arXiv preprint arXiv:1705.07957}, 2017.

\bibitem{ehn1}
H.~W. Engl, M.~Hanke, and A.~Neubauer.
\newblock {\em Regularization of Inverse Problems}.
\newblock Kluwer, Dordrecht, 1996.

\bibitem{erdogdu2015convergence}
Murat~A. Erdogdu and Andrea Montanari.
\newblock Convergence rates of sub-sampled newton methods.
\newblock In {\em Advances in Neural Information Processing Systems 28}, pages
  3034--3042. 2015.

\bibitem{fichtner}
A.~Fichtner.
\newblock {\em Full Seismic Waveform Modeling and Inversion}.
\newblock Springer, 2011.

\bibitem{fletcher2013practical}
R.~Fletcher.
\newblock {\em Practical methods of optimization}.
\newblock John Wiley \& Sons, 2013.

\bibitem{fox2016fast}
Colin Fox and Richard~A Norton.
\newblock Fast sampling in a linear-gaussian inverse problem.
\newblock {\em SIAM/ASA Journal on Uncertainty Quantification},
  4(1):1191--1218, 2016.

\bibitem{frank1956algorithm}
Marguerite Frank and Philip Wolfe.
\newblock An algorithm for quadratic programming.
\newblock {\em Naval research logistics quarterly}, 3(1-2):95--110, 1956.

\bibitem{gaooscher}
H.~Gao, S.~Osher, and H.~Zhao.
\newblock Quantitative photoacoustic tomography.
\newblock In {\em Mathematical Modeling in Biomedical Imaging II}, pages
  131--158. Springer, 2012.

\bibitem{ge2015escaping}
Rong Ge, Furong Huang, Chi Jin, and Yang Yuan.
\newblock Escaping from saddle points-online stochastic gradient for tensor
  decomposition.
\newblock In {\em COLT}, pages 797--842, 2015.

\bibitem{haber2004quasi}
E~Haber.
\newblock Quasi-newton methods for large-scale electromagnetic inverse
  problems.
\newblock {\em Inverse problems}, 21(1):305, 2004.

\bibitem{haasol}
E.~Haber, U.~Ascher, and D.~Oldenburg.
\newblock Inversion of 3{D} electromagnetic data in frequency and time domain
  using an inexact all-at-once approach.
\newblock {\em Geophysics}, 69:1216--1228, 2004.

\bibitem{HaberChungHermann2010}
E.~Haber, M.~Chung, and F.~Herrmann.
\newblock An effective method for parameter estimation with {PDE} constraints
  with multiple right-hand sides.
\newblock {\em SIAM J. Optimization}, 22:739--757, 2012.

\bibitem{haheas}
E.~Haber, S.~Heldmann, and U.~Ascher.
\newblock Adaptive finite volume method for distributed non-smooth parameter
  identification.
\newblock {\em Inverse Problems}, 23:1659--1676, 2007.

\bibitem{haber2001preconditioned}
Eldad Haber and Uri~M Ascher.
\newblock Preconditioned all-at-once methods for large, sparse parameter
  estimation problems.
\newblock {\em Inverse Problems}, 17(6):1847, 2001.

\bibitem{haber2000optimization}
Eldad Haber, Uri~M. Ascher, and Doug Oldenburg.
\newblock On optimization techniques for solving nonlinear inverse problems.
\newblock {\em Inverse problems}, 16(5):1263, 2000.

\bibitem{haber2014simultaneous}
Eldad Haber and Mathias Chung.
\newblock Simultaneous source for non-uniform data variance and missing data.
\newblock {\em arXiv preprint arXiv:1404.5254}, 2014.

\bibitem{haber2012effective}
Eldad Haber, Matthias Chung, and Felix Herrmann.
\newblock An effective method for parameter estimation with {PDE} constraints
  with multiple right-hand sides.
\newblock {\em SIAM Journal on Optimization}, 22(3):739--757, 2012.

\bibitem{hadamard1902sur}
Jacques Hadamard.
\newblock Sur les probl{\`e}mes aux d{\'e}riv{\'e}es partielles et leur
  signification physique.
\newblock {\em Princeton University Bulletin}, pages 49 -- 52, 1902.

\bibitem{hanke1}
M.~Hanke.
\newblock Regularizing properties of a truncated {N}ewton-cg algorithm for
  nonlinear inverse problems.
\newblock {\em Numer. Funct. Anal. Optim.}, 18:971--993, 1997.

\bibitem{hansen1998}
P.~C. Hansen.
\newblock {\em {Rank-Deficient and Discrete Ill-Posed Problems}}.
\newblock SIAM, 1998.

\bibitem{herman2009fundamentals}
Gabor~T Herman.
\newblock {\em Fundamentals of computerized tomography: image reconstruction
  from projections}.
\newblock Springer Science \& Business Media, 2009.

\bibitem{hel}
F.~Herrmann, Y.~Erlangga, and T.~Lin.
\newblock Compressive simultaneous full-waveform simulation.
\newblock {\em Geophysics}, 74:A35, 2009.

\bibitem{ito1990augmented}
Kazufumi Ito and Karl Kunisch.
\newblock The augmented lagrangian method for parameter estimation in elliptic
  systems.
\newblock {\em SIAM Journal on Control and Optimization}, 28(1):113--136, 1990.

\bibitem{jaggi2013revisiting}
Martin Jaggi.
\newblock Revisiting {F}rank-{W}olfe: Projection-free sparse convex
  optimization.
\newblock In {\em Proceedings of the 30th International Conference on Machine
  Learning (ICML-13)}, pages 427--435, 2013.

\bibitem{jain2013low}
Prateek Jain, Praneeth Netrapalli, and Sujay Sanghavi.
\newblock Low-rank matrix completion using alternating minimization.
\newblock In {\em Proceedings of the forty-fifth annual ACM symposium on Theory
  of computing}, pages 665--674. ACM, 2013.

\bibitem{jin2017escape}
Chi Jin, Rong Ge, Praneeth Netrapalli, Sham~M Kakade, and Michael~I Jordan.
\newblock How to escape saddle points efficiently.
\newblock {\em arXiv preprint arXiv:1703.00887}, 2017.

\bibitem{johnson2013accelerating}
Rie Johnson and Tong Zhang.
\newblock Accelerating stochastic gradient descent using predictive variance
  reduction.
\newblock In {\em Advances in Neural Information Processing Systems}, pages
  315--323, 2013.

\bibitem{kingma2014adam}
Diederik Kingma and Jimmy Ba.
\newblock Adam: A method for stochastic optimization.
\newblock {\em arXiv preprint arXiv:1412.6980}, 2014.

\bibitem{kumar2014GEOPemc}
R.~Kumar, C.~Da Silva, O.~Akalin, A.~Y. Aravkin, H.~Mansour, B.~Recht, and
  F.~J. Herrmann.
\newblock Efficient matrix completion for seismic data reconstruction.
\newblock Submitted to Geophysics on August 8, 2014., 08 2014.

\bibitem{lan2017conditional}
Guanghui Lan, Sebastian Pokutta, Yi~Zhou, and Daniel Zink.
\newblock Conditional accelerated lazy stochastic gradient descent.
\newblock In {\em ICML}. PMLR, 2017.

\bibitem{lan2016conditional}
Guanghui Lan and Yi~Zhou.
\newblock Conditional gradient sliding for convex optimization.
\newblock {\em SIAM Journal on Optimization}, 26(2):1379--1409, 2016.

\bibitem{levy2016power}
Kfir~Y Levy.
\newblock {The Power of Normalization: Faster Evasion of Saddle Points}.
\newblock {\em arXiv preprint arXiv:1611.04831}, 2016.

\bibitem{li2015accelerated}
Huan Li and Zhouchen Lin.
\newblock {Accelerated proximal gradient methods for nonconvex programming}.
\newblock In {\em {Advances in neural information processing systems}}, pages
  379--387, 2015.

\bibitem{liu1989limited}
Dong~C Liu and Jorge Nocedal.
\newblock On the limited memory {BFGS} method for large scale optimization.
\newblock {\em Mathematical programming}, 45(1-3):503--528, 1989.

\bibitem{liu2017noisy}
Mingrui Liu and Tianbao Yang.
\newblock {On Noisy Negative Curvature Descent: Competing with Gradient Descent
  for Faster Non-convex Optimization}.
\newblock {\em arXiv preprint arXiv:1709.08571}, 2017.

\bibitem{liu2017approximate}
Wenqing Liu, Jinglai Li, and Youssef~M Marzouk.
\newblock An approximate empirical bayesian method for large-scale
  linear-gaussian inverse problems.
\newblock {\em arXiv preprint arXiv:1705.07646}, 2017.

\bibitem{louis1992medical}
AK. Louis.
\newblock Medical imaging: state of the art and future development.
\newblock {\em Inverse Problems}, 8(5):709, 1992.

\bibitem{mazumder2011sparsenet}
Rahul Mazumder, Jerome~H Friedman, and Trevor Hastie.
\newblock Sparsenet: Coordinate descent with nonconvex penalties.
\newblock {\em Journal of the American Statistical Association},
  106(495):1125--1138, 2011.

\bibitem{menke2012geophysical}
W.~Menke.
\newblock {\em Geophysical data analysis: discrete inverse theory}.
\newblock Academic press, 2012.

\bibitem{mutny2016stochastic}
Mojm\`{i}r Mutn\`{y}.
\newblock Stochastic {S}econd-{O}rder {O}ptimization via von {N}eumann
  {S}eries.
\newblock {\em arXiv preprint arXiv:1612.04694}, 2016.

\bibitem{mutny2017parallel}
Mojm{\'\i}r Mutn{\`y} and Peter Richt{\'a}rik.
\newblock {Parallel Stochastic Newton Method}.
\newblock {\em arXiv preprint arXiv:1705.02005}, 2017.

\bibitem{natterer2001mathematical}
Frank Natterer and Frank W{\"u}bbeling.
\newblock {\em Mathematical methods in image reconstruction}.
\newblock SIAM, 2001.

\bibitem{na}
G.~A. Newman and D.~L. Alumbaugh.
\newblock Frequency-domain modelling of airborne electromagnetic responses
  using staggered finite differences.
\newblock {\em Geophys. Prospecting}, 43:1021--1042, 1995.

\bibitem{nocedal1980updating}
Jorge Nocedal.
\newblock Updating quasi-{N}ewton matrices with limited storage.
\newblock {\em Mathematics of computation}, 35(151):773--782, 1980.

\bibitem{nocedal2006numerical}
Jorge Nocedal and Stephen Wright.
\newblock {\em Numerical optimization}.
\newblock Springer Science \& Business Media, 2006.

\bibitem{olhash}
D.~Oldenburg, E.~Haber, and R.~Shekhtman.
\newblock 3{D} inverseion of multi-source time domain electromagnetic data.
\newblock {\em J. Geophysics}, 2013.
\newblock To appear.

\bibitem{of}
S.~Osher and R.~Fedkiw.
\newblock {\em Level Set Methods and Dynamic Implicit Surfaces}.
\newblock Springer, 2003.

\bibitem{osse}
S.~Osher and J.~Sethian.
\newblock Fronts propagating with curvature dependent speed: algorithms based
  on {H}amilton-{J}acobi formulations.
\newblock {\em J. Comp. Phys.}, 79:12--49, 1988.

\bibitem{pihakn}
A.~Pidlisecky, E.~Haber, and R.~Knight.
\newblock {RESINVM3D: A MATLAB 3D Resistivity Inversion Package}.
\newblock {\em Geophysics}, 72(2):H1--H10, 2007.

\bibitem{pilanci2015newton}
Mert Pilanci and Martin~J. Wainwright.
\newblock Newton sketch: A linear-time optimization algorithm with
  linear-quadratic convergence.
\newblock {\em arXiv preprint arXiv:1505.02250}, 2015.

\bibitem{reddi2016stochastic}
Sashank~J Reddi, Ahmed Hefny, Suvrit Sra, Barnabas Poczos, and Alex Smola.
\newblock Stochastic variance reduction for nonconvex optimization.
\newblock {\em arXiv preprint arXiv:1603.06160}, 2016.

\bibitem{rieder2005}
A.~Rieder.
\newblock Inexact {N}ewton regularization using conjugate gradients as inner
  iteration.
\newblock {\em SIAM J. Numer. Anal.}, 43:604--622, 2005.

\bibitem{rieder2010}
A.~Rieder and A.~Lechleiter.
\newblock Towards a general convergence theory for inexact {N}ewton
  regularizations.
\newblock {\em Numer. Math.}, 114(3):521--548, 2010.

\bibitem{rnkkda}
J.~Rohmberg, R.~Neelamani, C.~Krohn, J.~Krebs, M.~Deffenbaugh, and J.~Anderson.
\newblock Efficient seismic forward modeling and acquisition using simultaneous
  random sources and sparsity.
\newblock {\em Geophysics}, 75(6):WB15--WB27, 2010.

\bibitem{roosta2015randomized}
Farbod Roosta-Khorasani.
\newblock {\em Randomized algorithms for solving large scale nonlinear least
  squares problems}.
\newblock PhD thesis, University of British Columbia, 2015.

\bibitem{romassn1}
Farbod Roosta-Khorasani and Michael~W. Mahoney.
\newblock Sub-sampled {N}ewton methods {I}: Globally convergent algorithms.
\newblock {\em arXiv preprint arXiv:1601.04737}, 2016.

\bibitem{romassn2}
Farbod Roosta-Khorasani and Michael~W. Mahoney.
\newblock Sub-sampled {N}ewton methods {II}: Local convergence rates.
\newblock {\em arXiv preprint arXiv:1601.04738}, 2016.

\bibitem{roszas}
Farbod Roosta-Khorasani, G\'{a}bor~J. Sz\'{e}kely, and Uri Ascher.
\newblock Assessing stochastic algorithms for large scale nonlinear least
  squares problems using extremal probabilities of linear combinations of gamma
  random variables.
\newblock {\em SIAM/ASA Journal on Uncertainty Quantification}, 3(1):61--90,
  2015.

\bibitem{rodoas2}
Farbod Roosta-Khorasani, Kees van~den Doel, and Uri Ascher.
\newblock Data completion and stochastic algorithms for {PDE} inversion
  problems with many measurements.
\newblock {\em Electronic Transactions on Numerical Analysis}, 42:177--196,
  2014.

\bibitem{rodoas1}
Farbod Roosta-Khorasani, Kees van~den Doel, and Uri Ascher.
\newblock Stochastic algorithms for inverse problems involving {PDE}s and many
  measurements.
\newblock {\em SIAM J. Scientific Computing}, 36(5):S3--S22, 2014.

\bibitem{russell1988introduction}
B.~H. Russell.
\newblock {\em Introduction to seismic inversion methods}, volume~2.
\newblock Society of Exploration Geophysicists, 1988.

\bibitem{scharf1991statistical}
L.~L. Scharf.
\newblock {\em Statistical signal processing}, volume~98.
\newblock Addison-Wesley Reading, MA, 1991.

\bibitem{shalev2013stochastic}
Shai Shalev-Shwartz and Tong Zhang.
\newblock Stochastic dual coordinate ascent methods for regularized loss.
\newblock {\em The Journal of Machine Learning Research}, 14(1):567--599, 2013.

\bibitem{smvoz}
N.~C. Smith and K.~Vozoff.
\newblock Two dimensional {DC} resistivity inversion for dipole dipole data.
\newblock {\em IEEE Trans. on geoscience and remote sensing}, GE 22:21--28,
  1984.

\bibitem{sun2006optimization}
Wenyu Sun and Ya-Xiang Yuan.
\newblock {\em Optimization theory and methods: nonlinear programming},
  volume~1.
\newblock Springer Science \& Business Media, 2006.

\bibitem{sutskever2013importance}
Ilya Sutskever, James Martens, George Dahl, and Geoffrey Hinton.
\newblock On the importance of initialization and momentum in deep learning.
\newblock In {\em International conference on machine learning}, pages
  1139--1147, 2013.

\bibitem{tali}
X.-C. Tai and H.~Li.
\newblock A piecewise constant level set method for elliptic inverse problems.
\newblock {\em Appl. Numer. Math.}, 57:686--696, 2007.

\bibitem{tarantola2005inverse}
Albert Tarantola.
\newblock {\em Inverse problem theory and methods for model parameter
  estimation}.
\newblock SIAM, 2005.

\bibitem{tijmen2012rmsprop}
Tijmen Tieleman and Geoffrey Hinton.
\newblock Lecture 6.5-rmsprop: Divide the gradient by a running average of its
  recent magnitude.
\newblock {\em COURSERA: Neural Networks for Machine Learning}, 4, 2012.

\bibitem{tripuraneni2017stochastic}
Nilesh Tripuraneni, Mitchell Stern, Chi Jin, Jeffrey Regier, and Michael~I
  Jordan.
\newblock {Stochastic Cubic Regularization for Fast Nonconvex Optimization}.
\newblock {\em arXiv preprint arXiv:1711.02838}, 2017.

\bibitem{doas}
K.~van~den Doel and U.~M. Ascher.
\newblock On level set regularization for highly ill-posed distributed
  parameter estimation problems.
\newblock {\em J. Comp. Phys.}, 216:707--723, 2006.

\bibitem{doas1}
K.~van~den Doel and U.~M. Ascher.
\newblock Dynamic level set regularization for large distributed parameter
  estimation problems.
\newblock {\em Inverse Problems}, 23:1271--1288, 2007.

\bibitem{doas5}
K.~van~den Doel and U.~M. Ascher.
\newblock Dynamic regularization, level set shape optimization, and computed
  myography.
\newblock {\em Control and Optimization with Differential-Algebraic
  Constraints}, 23:315, 2012.

\bibitem{van2013lost}
Kees Van Den~Doel, Uri Ascher, and Eldad Haber.
\newblock The lost honour of $\ell_2$-based regularization.
\newblock {\em Radon Series in Computational and Applied Math}, 2013.

\bibitem{vogelbook}
C.~Vogel.
\newblock {\em Computational methods for inverse problem}.
\newblock SIAM, Philadelphia, 2002.

\bibitem{rundell1997inverse}
w.~Rundell and H.~W. Engl.
\newblock {\em Inverse problems in medical imaging and nondestructive testing}.
\newblock Springer-Verlag New York, Inc., 1997.

\bibitem{wang2015subsampled}
Chien-Chih Wang, Chun-Heng Huang, and Chih-Jen Lin.
\newblock Subsampled {H}essian {N}ewton methods for supervised learning.
\newblock {\em Neural computation}, 2015.

\bibitem{xu2017newton}
Peng Xu, Farbod Roosta-Khorasani, and Michael~W Mahoney.
\newblock Newton-type methods for non-convex optimization under inexact hessian
  information.
\newblock {\em arXiv preprint arXiv:1708.07164}, 2017.

\bibitem{xu2017second}
Peng Xu, Farbod Roosta-Khorasani, and Michael~W Mahoney.
\newblock Second-order optimization for non-convex machine learning: An
  empirical study.
\newblock {\em arXiv preprint arXiv:1708.07827}, 2017.

\bibitem{pyrrm_ssn_nonuni}
Peng Xu, Jiyan Yang, Farbod Roosta-Khorasani, Christopher R\'{e}, and
  Michael~W. Mahoney.
\newblock Sub-{S}ampled {N}ewton {M}ethods with {N}on-{U}niform {S}ampling.
\newblock In {\em Advances In Neural Information Processing Systems (NIPS)},
  pages 2530--2538, 2016.

\bibitem{ye2016revisiting}
Haishan Ye, Luo Luo, and Zhihua Zhang.
\newblock Revisiting sub-sampled newton methods.
\newblock {\em arXiv preprint arXiv:1608.02875}, 2016.

\bibitem{yuan}
Z.~Yuan and H.~Jiang.
\newblock Quantitative photoacoustic tomography: Recovery of optical absorption
  coefficient maps of heterogeneous media.
\newblock {\em Applied physics letters}, 88(23):231101--231101, 2006.

\bibitem{zeiler2012adadelta}
Matthew~D Zeiler.
\newblock Adadelta: an adaptive learning rate method.
\newblock {\em arXiv preprint arXiv:1212.5701}, 2012.

\end{thebibliography}
\end{document}